\documentstyle[11pt,twoside,a4wide]{article}

\makeatletter
\@addtoreset{equation}{section}
\makeatother

\def\parref#1{(\ref{#1})}

\def\Frac#1#2{\frac
{
 {\raise.6ex
 \hbox{$\displaystyle#1$}}
}
{
 {\lower.6ex
 \hbox{$\displaystyle#2$}}
 }
}

\def\morespace{\vspace{5pt}}
\def\sumslarge{{\displaystyle\sum_{s=0}^\infty}}
\def\begeq{\begin{equation}
\begin{array}{ll}
}
\def\endeq{\end{array}\end{equation}}

\input{psfig}
\def\intp{\int_0^\infty}

\def\bo{{\cal O}}
\def\Ai{{{\rm Ai}}}
\def\Bi{{{\rm Bi}}}
\def\sfrac#1#2{{{\lower.6ex
\hbox{$\scriptstyle#1$}}\over 
{\raise.7ex
\hbox{$\scriptstyle#2$}}}}

\begin{document}
 \title{
Numerical and Asymptotic Aspects of \\
Parabolic Cylinder Functions}

\author{Nico M. Temme\\
  CWI\\
  P.O. Box 94079\\
  1090 GB Amsterdam\\
  The Netherlands\\
e-mail: {\tt  nicot@cwi.nl}
}

 \maketitle
 \begin{abstract}
\noindent
Several uniform asymptotics expansions of the Weber
parabolic cylinder functions are considered, one group in terms of elementary
functions, another group in terms of Airy functions.  Starting point for the
discussion are  asymptotic expansions given earlier by F.W.J. Olver. Some of his
results are modified to improve the asymptotic properties and to enlarge the
intervals for using the expansions in numerical algorithms. 
Olver's results are obtained from the differential equation of the parabolic
cylinder functions; we mention how modified expansions can be obtained from integral
representations. Numerical tests are given for three expansions in terms of
elementary functions. In this paper only real values of the parameters will be
considered.
 \end{abstract}

\vskip 0.8cm \noindent
1991 Mathematics Subject Classification:
33C15, 41A60, 65D20.
\par\noindent
Keywords \& Phrases:
parabolic cylinder functions,
uniform asymptotic expansion,
Airy-type expansions, 
numerical evaluation of special functions.


\section{Introduction}%
The solutions of the differential equation
\begin{equation}
\frac{d^2y}{dz^2}-\left(\sfrac14z^2+a\right)y=0 \label{p1}
\end{equation}
are associated with the parabolic cylinder in harmonic analysis; see
\cite{weber}. The solutions are called parabolic cylinder functions
and are entire functions  of
$z$. Many properties in connection with physical applications are given in 
\cite{buch}.

As in \cite{abst}, Chapter 19, we denote 
two standard solutions of \parref{p1} by $U(a,z), V(a,z)$. These solutions 
are given by the representations

\begeq
U(a,z)&=\sqrt{{\pi}}2^{-\sfrac12a}\left[
\Frac{2^{-\sfrac14}\,y_1(a,z)}{\Gamma(\sfrac34+\sfrac12a)}
-
\Frac{2^{\sfrac14}\,y_2(a,z)}{\Gamma(\sfrac14+\sfrac12a)}\right],\\ 
V(a,z)&=\Frac{\sqrt{{\pi}}2^{-\sfrac12a}}{\Gamma(\sfrac12-a)}\left[
\tan\pi\left(\sfrac12a+\sfrac14\right)\Frac{2^{-\sfrac14}\,y_1(a,z)}
{\Gamma(\sfrac34+\sfrac12a)}
+
\cot\pi\left(\sfrac12a+\sfrac14\right)\Frac{2^{\sfrac14}\,y_2(a,z)}
{\Gamma(\sfrac14+\sfrac12a)}\right],\label{p2}
\endeq

where
\begeq
y_1(a,z)&= e^{\sfrac14z^2}{}_1F_1\left(-\sfrac12a+\sfrac14,\sfrac12;-\sfrac12z^2\right)
=e^{-\sfrac14z^2}{}_1F_1\left(\sfrac12a+\sfrac14,\sfrac12;\sfrac12z^2\right), \\
y_2(a,z)&= ze^{-\sfrac14z^2}{}_1F_1\left(\sfrac12a+\sfrac34,\sfrac32;\sfrac12z^2\right)
=ze^{\sfrac14z^2}{}_1F_1\left(-\sfrac12a+\sfrac34,\sfrac32;-\sfrac12z^2\right)\label{p3}
\endeq

and the confluent hypergeometric function is defined by
\begin{equation}
{}_1F_1(a,c;z)=\sum_{n=0}^\infty \frac{(a)_n}{(c)_n}\,\frac{z^n}{n!},\label{p4}
\end{equation}
with $(a)_n=\Gamma(a+n)/\Gamma(a), n=0,1,2,\ldots$. 

Another notation found in the literature is
$$D_\nu(z)=U(-\nu-\sfrac12,z).$$
There is a relation with the Hermite polynomials. We have
\begeq
U\left(-n-\sfrac12,z\right)&=2^{-n/2}e^{-\sfrac14z^2}H_n(z/\sqrt{{2}}),\\
V\left( n+\sfrac12,z\right)&=2^{-n/2}e^{\sfrac14z^2}(-i)^nH_n(iz/\sqrt{{2}}).\label{p5}
\endeq
Other special cases are error functions and Fresnel integrals.

The Wronskian relation between $U(a,z)$ and $V(a,z)$ reads:
\begin{equation}
U(a,z)V'(a,z)-U'(a,z)V(a,z)=\sqrt{{2/\pi}},\label{p6}
\end{equation}
which shows that $U(a,z)$ and $V(a,z)$ are independent solutions of \parref{p1}
for all values of  $a$. Other relations are
\begeq
U(a,z)&=\Frac{\pi}{\cos^2\pi a\,\Gamma(\sfrac12+a)}
\left[V(a,-z)-\sin\pi a\,V(a,z)\right],\\
V(a,z)&=\Frac{\Gamma(\sfrac12+a)}\pi
\left[\sin\pi a\,U(a,z)+U(a,-z)\right].\label{p7}
\endeq
The functions $y_1(a,z)$ and $y_2(a,z)$ are the simplest even and odd solutions
of \parref{p1} and the Wronskian of this pair equals 1. From a numerical point of
view, the pair $\{y_1,y_2\}$ is not a satisfactory pair \cite{mil1}, 
because they have almost the same asymptotic behavior at infinity.
The behavior of $U(a,z)$ and  $V(a,z)$ is, for large positive $z$ and $z\gg|a|$:
\begeq
U(a,z)&=e^{-\sfrac14z^2}z^{-a-\sfrac12}\left[1+\bo\left(z^{-2}\right)\right],\\
V(a,z)&=\sqrt{{2/\pi}}e^{\sfrac14z^2}z^{a-\sfrac12}\left[1+\bo\left(z^{-2}\right)\right].
\label{p8}
\endeq
Clearly, numerical computations of $U(a,z)$ that are based on the representations
in \parref{p2} should be done with great care, because of the loss of accuracy if
$z$ becomes large.

Equation \parref{p1} has two turning points at $\pm 2\sqrt{{-a}}$. For real parameters 
they become important if $a$ is negative, and the asymptotic behavior of the 
solutions of \parref{p1} as $a\to-\infty$  changes significantly if $z$ crosses the
turning points. At these points Airy functions are needed. By changing the
parameters it is not difficult to verify that $U(-\sfrac12\mu^2,\mu t\sqrt{2})$
and $V(-\sfrac12\mu^2,\mu t\sqrt{2})$ satisfy the simple equation
\begin{equation}
\frac{d^2y}{dt^2}-\mu^4\left(t^2-1\right)y=0, \label{p9}
\end{equation}
with turning points at $t=\pm1$. For physical applications, negative
$a-$values are most important (with special case the real Hermite polynomials,
see \parref{p5}). For positive $a$ we can use the notation
$U(\sfrac12\mu^2,\mu t\sqrt{2})$ and $V(\sfrac12\mu^2,\mu t\sqrt{2})$, which satisfy the
equation
\begin{equation}
\frac{d^2y}{dt^2}-\mu^4\left(t^2+1\right)y=0. \label{p10}
\end{equation}

The purpose of this paper is to give several asymptotic expansions of 
$U(a,z)$ and  $V(a,z)$ that can be used for computing these functions for the
case that at least one of the real parameters is large. In \cite{olpar}
an extensive collection of asymptotic expansions for the parabolic cylinder
functions as $|a|\to\infty$ has been derived from the differential equation
\parref{p1}. The expansions are valid for complex values of the parameters and are
given in terms of elementary functions and Airy functions.
In Section 2 we mention several expansions in terms of elementary functions derived
by Olver and modify some his results in order to improve the asymptotic properties of
the expansions, to enlarge the intervals for using the expansions in numerical
algorithms, and to get new recursion relations for the coefficients of the
expansions. In Section 3 we give similar results for expansions in terms of Airy
functions. In Section 4 we give information on how to obtain the modified results by
using integral representations of the parabolic cylinder functions. 
Finally we give numerical tests for three expansions in terms of
elementary functions, with a few number of terms in the expansions. Only real
parameters are considered in this paper.
\subsection{Recent literature on numerical algorithms}%
Recent papers on numerical algorithms for the parabolic cylinder 
functions are given in \cite{taub} (Fortran; $U(n,x)$ for natural 
$n$ and positive $x$) and \cite{seggil} (Fortran; $U(a,x), V(a,x)$, 
$a$ integer and half-integer and $x\ge 0$). The methods are based on backward 
and forward recursion.  

\cite{baker} gives programs in ${{\cal C}}$ 
for $U(a,x), V(a,x)$, and uses representations in terms of 
the confluent hypergeometric functions and 
asymptotic expressions, including those involving Airy functions.
\cite{zhjin} gives Fortran programs for computing
$U(a,z), V(a,z)$ with real orders and real argument, 
and for half-integer order and
complex argument, The methods are based on recursions, Maclaurin series and
asymptotic expansions. They refer also to \cite{blannew}
for the evaluation of $U(-ia,ze^{\sfrac14\pi i})$ for real $a$ and $z$ 
(this function is a solution of the differential 
equation $y''+(\sfrac14z^2-a)y=0$).
\cite{tomp} uses series expansions and numerical quadrature; 
Fortran and
${{\cal C}}$ programs are given, and \emph{Mathematica} cells to make
graphical and numerical objects.

\emph{Maple}  \cite{maple} has algorithms for hypergeometric functions, which can be used in
\parref{p2} and \parref{p3}. 
\emph{Mathematica} \cite{math} refers for the parabolic cylinder functions
to their programs for the hypergeometric functions
 and the same advice is given 
in \cite{numrec}.

\section{Expansions in terms of elementary functions}%
\subsection{The case $\quad a\le0,\quad z>2\sqrt{{-a}},\quad -a+z\gg0$.}
Olver's expansions in
terms of elementary functions are all based on the expansion 
O-(4.3)
\footnote{We refer to Olver's equations by writing O-(4.3), and so on.}
\begin{equation}
U\left(-\sfrac12\mu^2,\mu t\sqrt{2}\right)\sim 
\frac{g(\mu)\,e^{-\mu^2\xi}}{(t^2-1)^{\frac14}}
\sum_{s=0}^\infty\frac{{\cal A}_s(t)}{\mu^{2s}},\label{o1}
\end{equation}
as $\mu\to\infty$, uniformly with respect to  $t\in[1+\varepsilon,\infty)$;
$\varepsilon$ is a small positive number and $\xi$ is given by
\begin{equation} \xi=\sfrac12t\sqrt{{t^2-1}}-\sfrac12\ln\left[t+\sqrt{{t^2-1}}\right].\label{o2}
\end{equation}
The expansion is valid for complex parameters in large domains of the $\mu-$ and
$t-$planes; details on these domains are not given here.

The coefficients ${\cal A}_s(t)$ are given by the recursion relation
\begin{equation}
{\cal A}_{s+1}(t)=\sfrac12\,\frac1{\sqrt{{t^2-1}}}\,\frac{d{\cal
A}_s(t)}{dt}+\sfrac18\int_{c_{s+1}}^t
\frac{3u^2+2}{(u^2-1)^\frac52}\,{\cal A}_s(u)\,du,\quad {\cal A}_0(t)=1,\label{o3}
\end{equation}
where the constants $c_s$ can be chosen in combination with the choice of
$g(\mu)$. Olver chose the constants such that
\begin{equation}
{\cal A}_s(t)=\frac{u_s(t)}{(t^2-1)^{3s/2}},\label{o4}
\end{equation}
where the $u_s(t)$ are polynomials in $t$ of degree $3s$, ($s$ odd), $3s-2$ ($s$
even, $s\ge2$). The first few are
$$u_0(t)=1,\quad u_1(t)=\frac{t(t^2-6)}{24}, \quad
u_2(t)=\frac{-9t^4+249t^2+145}{1152},$$  
and they satisfy the recurrence relation
\begin{equation}
(t^2-1)u'_s(t)-3stu_s(t)=r_{s-1}(t),\label{o5}
\end{equation}
where
$$8r_s(t)=(3t^2+2)u_s(t)-12(s+1)tr_{s-1}(t)+4(t^2-1)r'_{s-1}(t).$$
The quantity $g(\mu)$ in \parref{o1} is only available in the form of an asymptotic
expansion:
\begin{equation}
g(\mu)\sim h(\mu)
\left(\sum_{s=0}^\infty\frac{g_s}{\mu^{2s}}\right)^{-1},\label{o6}
\end{equation}
where
\begin{equation}
h(\mu)=2^{-\sfrac14\mu^2-\sfrac14}e^{-\sfrac14\mu^2}
\mu^{\sfrac12\mu^2-\sfrac12},\label{o7}
\end{equation}
$$g_0=1,\quad g_1=\sfrac1{24},\quad g_3=-\sfrac{2021}{207360},\quad g_{2s}=0\ 
(s=1,2,\ldots),$$
and in general
\begin{equation}
g_s=\lim_{t\to\infty}{\cal A}_s(t).\label{o8}
\end{equation}
\subsubsection{Modified expansions}%
We modify the expansion in \parref{o1} by writing
\begin{equation}
U\left(-\sfrac12\mu^2,\mu t\sqrt{2}\right)=
\frac{h(\mu)\,e^{-\mu^2\xi}}{(t^2-1)^{\frac14}}F_\mu(t),\quad
F_\mu(t)\sim\sum_{s=0}^\infty\frac{\phi_s(\tau)}{\mu^{2s}},\label{o9}
\end{equation}
where $h(\mu)$ and $\xi$ are as before, and
\begin{equation}
\tau=\sfrac12\left[\frac{t}{\sqrt{{t^2-1}}}-1\right].\label{o10}
\end{equation}
The coefficients $\phi_s(\tau)$ are polynomials in $\tau$, $\phi_0(\tau)=1$, and are
given by the recursion
\begin{equation}
\phi_{s+1}(\tau)=-4\tau^2(\tau+1)^2\frac{d}{d\tau}\phi_s(\tau)
-\sfrac14\int_0^\tau\left(20u^2+20u+3\right)\phi_{s}(u)\,du.
\label{o11}
\end{equation}
This recursion follows from \parref{o3} by substituting
$t=(\tau+\sfrac12)/\sqrt{{\tau(\tau+1)}}$, which is the inverse of the relation in
\parref{o10}. 
Explicitly,

\begeq
\phi_0(\tau)&=   1,\morespace\\
\phi_1(\tau)&=   -\Frac{\tau}{12}  (20\tau^2+30\tau+9),\morespace\\
\phi_2(\tau)&=    \Frac{\tau^2}{288}(6160\tau^4+18480\tau^3+19404\tau^2+8028\tau+945),\morespace\\
\phi_3(\tau)&=   -\Frac{\tau^3}{51840}(27227200\tau^6+ 122522400\tau^5 + 220540320\tau^4+\morespace\\
       &   \quad\quad   200166120\tau^3+ 94064328\tau^2+
          20545650\tau+1403325),\label{o12}
\endeq
where $\tau$ is given in \parref{o10}. Observe that $\lim_{t\to\infty}\tau(t)=0$
and that all shown coefficients $\phi_s(\tau)$ vanish at infinity for $s>0$.
These properties of $\phi_s(\tau)$ follow by taking different constants $c_s$ 
than Olver did in \parref{o3}. In fact we have the relation
$$\sum_{s=0}^\infty
\frac{g_{s}}{\mu^{2s}} \sum_{s=0}^\infty\frac{\phi_s(\tau)}{\mu^{2s}}\sim
\sum_{s=0}^\infty\,\frac{{u}_s(t)}
{(t^2-1)^{\frac32s}\mu^{2s}},$$
where the first series appears in \parref{o6}. Explicitly,
\begin{equation}
u_s(t)=\left(t^2-1\right)^{\sfrac32s}\,\sum_{j=0}^s g_{s-j}\phi_j(\tau).\label{o13}
\end{equation}
The relation \parref{o13} can easily be verified for the early coefficients, but
it holds because of the unicity of Poincar\'e-type asymptotic expansions.

The expansion in \parref{o9} has several advantages compared with \parref{o1}.
\begin{itemize}
\item
In the recursion relation \parref{o5}, both $u_s$ and $u'_s$ occur in the
left-hand side. By using computer algebra it is not difficult to compute any
number of coefficients $u_s$, but the relation for the polynomials $\phi_s(\tau)$
is simpler than  \parref{o5}, with this respect.
\item
The quantity $h(\mu)$ in \parref{o9} is defined as an
exact relation, and not, as $g(\mu)$ in \parref{o1}, by an asymptotic expansion
(cf. \parref{o6}).
\item
Most important, the expansion in \parref{o9} has a double asymptotic property: it
holds if one or both parameters $t$ and  $\mu$ are large, and not only if
$\mu$ is large.
\end{itemize}

\medskip

For the function $V(a,z)$ we have 
\begin{equation}
V\left(-\sfrac12\mu^2,\mu t\sqrt{2}\right)=
 \frac{       e^{\mu^2\xi}}{\mu\,\sqrt{{\pi}}\,h(\mu)(t^2-1)^{\frac14}}
\,P_\mu(t),\quad
P_\mu(t)\sim\sum_{s=0}^\infty(-1)^s\frac{\phi_s(\tau)}{\mu^{2s}},\label{o14}
\end{equation}
where the $\phi_s(\tau)$ are the same as in \parref{o9}. This expansion is a
modification of O-(11.19) (see also O-(2.12)).

For the
derivatives we can use the identities
\begeq
\Frac d{dt}\Frac{e^{-\mu^2\xi}}{(t^2-1)^{\frac14}}
F_\mu(t)
&=-\mu^2(t^2-1)^{\sfrac14}
e^{-\mu^2\xi}\,G_\mu(t),\quad G_\mu(t)\sim
\sumslarge \Frac{\psi_s(\tau)}{\mu^{2s}},\\
\Frac d{dt}\Frac{e^{+\mu^2\xi}}{(t^2-1)^{\frac14}}
P_\mu(t)
&=+\mu^2(t^2-1)^{\sfrac14}
e^{+\mu^2\xi}\,Q_\mu(t),\quad Q_\mu(t)\sim
\sumslarge(-1)^s\Frac{\psi_s(\tau)}{\mu^{2s}}.\label{o15}
\endeq
The coefficients $\psi_s$ can be obtained from the relation
\begin{equation}
\psi_s(t)=\phi_s(\tau)+2\tau(\tau+1)(2\tau+1)\phi_{s-1}(\tau)+8\tau^2(\tau+1)^2\frac{d\phi_{s-
1}(\tau)}{d\tau},\label{o16}
\end{equation}
$s=0,1,2,\ldots\,$. The first few are
\begeq
\psi_0(t)&=   1,\morespace\\
\psi_1(t)&=   \Frac{\tau}{12}  (28\tau^2+42\tau+15),v\\
\psi_2(t)&=   -\Frac{\tau^2}{288}(7280\tau^4+21840\tau^3+23028\tau^2+9684\tau+1215),\morespace\\
\psi_3(t)&=\Frac{\tau^3}{51840}(30430400\tau^6+136936800\tau^5+246708000\tau^4+\morespace\\
        &   \quad\quad
224494200\tau^3+106122312\tau^2+23489190\tau+1658475).\label{o17}
\endeq

This gives the modifications (see O-(4.13))
\begin{equation}
U'\left(-\sfrac12\mu^2,\mu t\sqrt{2}\right)=
-\frac{\mu}{\sqrt{2}}h(\mu)(t^2-1)^{\sfrac14}e^{-\mu^2\xi}\,G_\mu(t),\quad G_\mu(t)\sim
\sum_{s=0}^\infty\frac{\psi_s(\tau)}{\mu^{2s}} \label{o18}
\end{equation}
and
\begin{equation}
V'\left(-\sfrac12\mu^2,\mu t\sqrt{2}\right)=
\frac{(t^2-1)^{\sfrac14}e^{\mu^2\xi}}{\sqrt{{2\pi}} h(\mu)}\,Q_\mu(t),
\quad Q_\mu(t)\sim
\sum_{s=0}^\infty(-1)^s\frac{\psi_s(\tau)}{\mu^{2s}}. \label{o19}
\end{equation}

\noindent
{\bf Remark 2.1.\ }
The functions $F_\mu(t), G_\mu(t), P_\mu(t)$ and $Q_\mu(t)$ introduced in
the asymptotic representations satisfy the following exact relation:
\begin{equation}
F_\mu(t)\,Q_\mu(t)+G_\mu(t)\,P_\mu(t)=2.\label{o20}
\end{equation}
This follows from the Wronskian relation (1.6).
The relation in
\parref{o20} provides a convenient possibility for checking the accuracy
in numerical algorithms that use the asymptotic expansions of
$F_\mu(t), G_\mu(t), P_\mu(t)$ and $Q_\mu(t)$. 

%
%
%
\subsection{The case $\quad a\le0,\quad z<-2\sqrt{{-a}},\quad -a-z\gg0$.}
For this case we mention the modification of 
O-(11.16). That is, for $t\ge1+\varepsilon$ we have the representations
\begeq
&U\left(-\sfrac12\mu^2,-\mu t\sqrt{2}\right)=
\Frac{h(\mu)}{(t^2-1)^{\frac14}}\quad\times\\
&\quad\quad\left[\sin\left(\sfrac12\pi\mu^2\right)e^{-\mu^2\xi}\,F_\mu(t)
+\Frac{\Gamma(\sfrac12+\sfrac12\mu^2)\cos\left(\sfrac12\pi\mu^2\right)}{
\mu\,\sqrt{{\pi}}\,h^2(\mu)}e^{\mu^2\xi}\,P_\mu(t)
\right],\label{o21}
\endeq 
where $F_\mu(t)$ and $P_\mu(t)$ have the expansions given in
\parref{o9} and \parref{o14}, respectively.
An expansion for $V(-\sfrac12\mu^2,-\mu t\sqrt{2})$ follows from the second line
in (1.7), \parref{o9} and \parref{o21}. A few manipulations give
\begeq
&V\left(-\sfrac12\mu^2,-\mu t\sqrt{2}\right)=
\Frac{h(\mu)}{(t^2-1)^{\frac14}\Gamma(\sfrac12+\sfrac12\mu^2)}\quad\times\\
&\quad\quad\left[\cos\left(\sfrac12\pi\mu^2\right)e^{-\mu^2\xi}\,F_\mu(t)
-\Frac{\Gamma(\sfrac12+\sfrac12\mu^2)\sin\left(\sfrac12\pi\mu^2\right)}{
\mu\,\sqrt{{\pi}}\,h^2(\mu)}e^{\mu^2\xi}\,P_\mu(t)
\right].\label{o22}
\endeq
Expansions for the derivatives follow from the identities in \parref{o15}. If
$a=-\sfrac12\mu^2=-n-\sfrac12, n = 0,1,2,\ldots$, the cosine in
\parref{o21}  vanishes, and, hence, the dominant part vanishes. This is the
Hermite case, cf. (1.5).

\subsection{The case $\quad  a\ll0,\quad -2\sqrt{{-a}}<z<2\sqrt{{-a}}$.}
For negative $a$ and $-1<t<1$ the expansions are essentially different, because
now oscillations with respect to $t$ occur. We have (O-(5.11) and O-(5.23))
\begeq&U\left(-\sfrac12\mu^2,\mu t\sqrt{2}\right)\sim
\Frac{2g(\mu)}{(1-t^2)^{\sfrac14}}\quad\times\\
&\quad\quad\left[\cos\left(\mu^2\eta-\sfrac14\pi\right)
\sumslarge\Frac{(-1)^su_{2s}(t)}{(1-t^2)^{3s}\mu^{4s}}
-\sin\left(\mu^2\eta-\sfrac14\pi\right)
\sumslarge\Frac{(-1)^su_
{2s+1}(t)}{(1-t^2)^{3s+\frac32}\mu^{4s+2}}\right],\label{o23}
\endeq
with $u_s(t)$ defined in \parref{o5} and $g(\mu)$ in \parref{o6}, and
\begeq&
U'\left(-\sfrac12\mu^2,\mu t\sqrt{2}\right)\sim
\mu\sqrt{2} g(\mu)(1-t^2)^{\sfrac14}\quad\times\\
&\quad\quad\left[\sin\left(\mu^2\eta-\sfrac14\pi\right)
\sumslarge\Frac{(-1)^sv_{2s}(t)}{(1-t^2)^{3s}\mu^{4s}}
+\cos\left(\mu^2\eta-\sfrac14\pi\right)
\sumslarge\Frac{(-1)^sv_
{2s+1}(t)}{(1-t^2)^{3s+\frac32}\mu^{4s+2}}\right],\label{o24}
\endeq
as $\mu\to\infty$, uniformly with respect to $|t|\le1-\varepsilon$, where the
coefficients $v_s$ are given by (see O-(4.15))
\begin{equation}
v_s(t)=u_s(t)+\sfrac12tu_{s-1}(t)-r_{s-2}(t), \label{o25}
\end{equation}
and
\begin{equation}
\eta=\sfrac12\arccos t-\sfrac12t\sqrt{{1-t^2}}.\label{o26}
\end{equation}

For the function $V(a,z)$ we have (O-(11.20) and O(2.12))
\begeq&V\left(-\sfrac12\mu^2,\mu t\sqrt{2}\right)\sim
\Frac{2g(\mu)}{\Gamma(\sfrac12+\sfrac12\mu^2)(1-t^2)^{\sfrac14}}\quad\times\\
&\quad\quad\left[\cos\left(\mu^2\eta+\sfrac14\pi\right)
\sumslarge\Frac{(-1)^su_{2s}(t)}{(1-t^2)^{3s}\mu^{4s}}
-\sin\left(\mu^2\eta+\sfrac14\pi\right)
\sumslarge\Frac{(-1)^su_
{2s+1}(t)}{(1-t^2)^{3s+\frac32}\mu^{4s+2}}\right],\\
&V'\left(-\sfrac12\mu^2,\mu t\sqrt{2}\right)\sim
\Frac{\mu\sqrt{2} g(\mu)(1-t^2)^{\sfrac14}}{\Gamma(\sfrac12+\sfrac12\mu^2)}\quad\times\\
&\quad\quad\left[\sin\left(\mu^2\eta+\sfrac14\pi\right)
\sumslarge\Frac{(-1)^sv_{2s}(t)}{(1-t^2)^{3s}\mu^{4s}}
+\cos\left(\mu^2\eta+\sfrac14\pi\right)
\sumslarge\Frac{(-1)^sv_
{2s+1}(t)}{(1-t^2)^{3s+\frac32}\mu^{4s+2}}\right],\label{o27}
\endeq

By using the Wronskian relation (1.6) it follows that we have the following
asymptotic identity
\begeq
&\sumslarge\Frac{(-1)^su_{2s}(t)}{(1-t^2)^{3s}\mu^{4s}}\ 
\sumslarge \Frac{(-1)^sv_{2s}(t)}{(1-t^2)^{3s}\mu^{4s}}
+\morespace\\
&\sumslarge\Frac{(-1)^su_{2s+1}(t)}{(1-t^2)^{3s+\sfrac32}\mu^{4s+2}}\ 
\sumslarge\Frac{(-1)^sv_{2s+1}(t)}{(1-t^2)^{3s+\sfrac32}\mu^{4s+2}}\morespace\\
&\sim\Frac{\Gamma(\sfrac12+\sfrac12\mu^2)}{2\mu\sqrt{{\pi}}g^2(\mu)}
\sim1-\Frac1{576\mu^4}+\Frac{2021}{2488320\mu^8}+\ldots\ .
\label{o28}
\endeq
\subsubsection{Modified expansions}%
We can give modified versions based on our earlier modifications, with 
$g(\mu)$ replaced with $h(\mu)$, and so on.  
Because in the present case $t$ belongs to a finite domain,
the modified expansions don't have the double asymptotic
property. We prefer Olver's versions for this case. 

This completes the description of $U(a,z), U'(a,z), V(a,z),
V'(a,z)$ in terms of elementary functions for negative values of $a$. 

\subsection{The case $\quad a\ge0, \quad z\ge0,\quad a+z\gg0$.}
For positive values of $a$
the asymptotic behavior is rather simple because no oscillations occur now.
The results follow from Olver's expansions
O-(11.10) and O-(11.12). The modified forms are
\begin{equation}
U\left(\sfrac12\mu^2,\mu t\sqrt{2}\right)=
\frac{\widetilde h(\mu)\,e^{-\mu^2{\widetilde\xi}}}
{(t^2+1)^{\frac14}}\,\widetilde F_\mu(t),\quad
\widetilde F_\mu(t)\sim\sum_{s=0}^\infty\,(-1)^s\,
\frac{\phi_s(\widetilde\tau)}{\mu^{2s}},\label{o29}
\end{equation}
where 
\begin{equation}
\widetilde\xi=\sfrac12\left[t\sqrt{{1+t^2}}+\ln\left(t+\sqrt{{1+t^2}}\right)\right],
\label{o30}
\end{equation}
\begin{equation}\widetilde h(\mu)=e^{\sfrac14\mu^2}\,\mu^{-\sfrac12\mu^2-\sfrac12}\,
2^{\sfrac14\mu^2-\sfrac14}.\label{o31}
\end{equation}
The coefficients $\phi_s$ in \parref{o29} are the same as in \parref{o9},
with $\tau$ replaced by
\begin{equation}
\widetilde\tau=\sfrac12\left[\frac{t}{\sqrt{{1+t^2}}}-1\right].\label{o32}
\end{equation}
For the derivative we have
\begin{equation}
U'\left(\sfrac12\mu^2,\mu t\sqrt{2}\right)=-\frac1{\sqrt{2}}\mu
\widetilde h(\mu)(1+t^2)^{\sfrac14}\,e^{-\mu^2{\widetilde\xi}}
{\widetilde G}_a(z),\quad
{\widetilde G}_a(z)\sim
\sum_{s=0}^\infty\,(-1)^s\,\frac{\psi_s(\widetilde\tau)}{\mu^{2s}},\label{o33}
\end{equation} 
where $\psi_s(\widetilde\tau)$ is given in \parref{o16}, with $\widetilde\tau$
defined in \parref{o32}.

\subsection{The case $\quad a\ge0, \quad z\le0,\quad a-z\gg0$.}
Olver's expansion O-(11.10) and O(11.12) cover both cases $z\ge0$ and $z\le0$. 
We have the modified expansions 
\begeq
U\left(\sfrac12\mu^2,-\mu t\sqrt{2}\right)
&=\Frac{\sqrt{{2\pi}}}{\Gamma(\sfrac12+\sfrac12\mu^2)}
\Frac{h(\mu)e^{\mu^2{\widetilde\xi}}}{(1+t^2)^{\frac14}}\,
\widetilde P_\mu(t),\\ 
U'\left(\sfrac12\mu^2,-\mu t\sqrt{2}\right)
&=-\Frac\mu{\sqrt{2}}\Frac{\sqrt{{2\pi}}}{\Gamma(\sfrac12+\sfrac12\mu^2)}
 h(\mu)e^{\mu^2{\widetilde\xi}}(1+t^2)^{\sfrac14}\,
\widetilde Q_\mu(t),\label{o34}
\endeq
where
$$
\widetilde P_\mu(t)\sim\sum_{s=0}^\infty \frac{\phi_s  (\widetilde\tau)}{\mu^{2s}},\quad
\widetilde Q_\mu(t)\sim\sum_{s=0}^\infty \frac{\psi_s(\widetilde\tau)}{\mu^{2s}}.
$$
In Section 4.1.2 we give details on the derivation of these expansions.

\noindent
{\bf Remark 2.2.\ }
By using the second relation in (1.7), the representations for $V(a,z)$ and
$V'(a,z)$ for positive $a$ can be obtained from the results for $U(a,z)$ and
$U'(a,z)$ in \parref{o29}, \parref{o33} and \parref{o34}.

\noindent
{\bf Remark 2.3.\ }
The functions $\widetilde F_\mu(t), \widetilde G_\mu(t), 
\widetilde P_\mu(t)$ and $\widetilde Q_\mu(t)$ introduced 
in \parref{o29}, \parref{o32} and \parref{o34} satisfy the following exact relation:
\begin{equation}
\widetilde F_\mu(t)\widetilde Q_\mu(t)+
\widetilde G_\mu(t)\widetilde P_\mu(t)=2.\label{o35}
\end{equation}
This follows from the Wronskian relation
$$U(a,z)U'(a,-z)+U'(a,z)U(a,-z)=-\frac{\sqrt{{2\pi}}}{\Gamma(a+\sfrac12)}.$$
See also Remark 2.1.

\noindent
{\bf Remark 2.4.\ }
The expansions of Sections 2.4 and 2.5 have the double asymptotic property: they
are valid if the $a+|z|\to\infty$. In Sections 2.4 and 2.5 we consider the cases
$z\ge0$ and
$z\le0$, respectively, as two separate cases. Olver's corresponding expansions
O-(11.10) and O-(11.12) cover both cases and are valid for $-\infty<t<\infty$.
As always, in Olver's expansions large values of $\mu$ are needed, 
whatever the size of $t$.

In Figure 1 we show the domains in the $t,a-$plane where the
various expansions of $U(a,z)$ of this section are valid.

\vspace*{0.3cm}

\centerline{\protect\hbox{\psfig{file=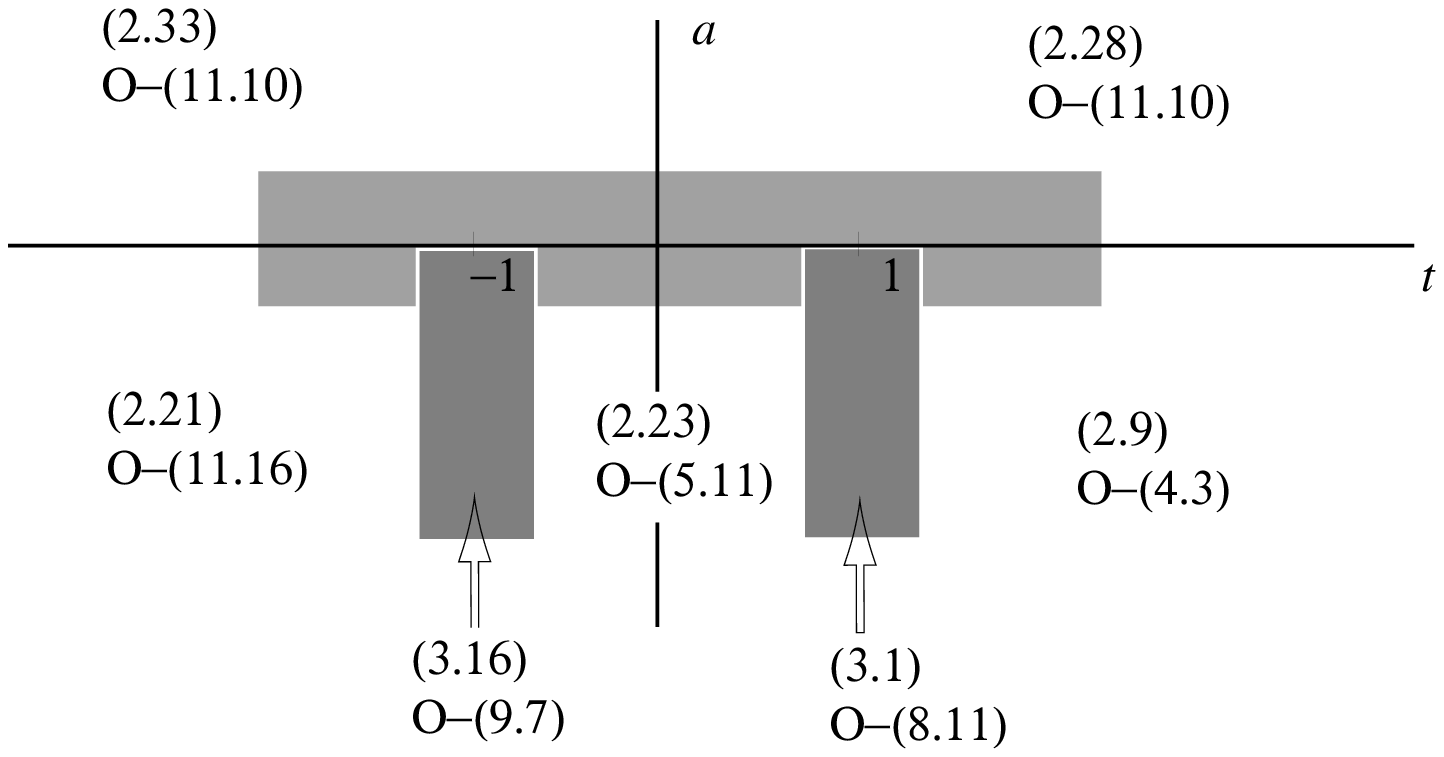,width=12cm}}}

\begin{quote}
{\bf Figure 1.}\quad
Regions for the modified asymptotic expansions of $U(a,z)$ given in 
Section 2 and the Airy-type expansions of Section 3 (which are valid in 
much larger domains than those indicated by the arrows).
\end{quote}

\section{Expansions in terms of Airy functions}%
The Airy-type expansions are needed if $z$ runs through an interval
containing one of the turning points $\pm2\sqrt{{-a}}$, that is,
$t=\pm1$. 
\noindent
\subsection{The case $\quad a\ll0,\quad z\ge0$.}
We summarize the basic results O-(8.11), O-(8.15) and O-(11.22) (see also 
O-(2.12)):
\begin{equation}
U\left(-\sfrac12\mu^2,\mu t\sqrt{2}\right)=
2\pi^{\sfrac12}\mu^{\frac13}g(\mu)\phi(\zeta)\left[
\Ai(\mu^{\frac43}\zeta)\,A_\mu(\zeta)+
\frac{\Ai'(\mu^{\frac43}\zeta)}{\mu^{\frac83}}\,B_\mu(\zeta)\right],\label{a1}
\end{equation} 
\begin{equation}
U'\left(-\sfrac12\mu^2,\mu t\sqrt{2}\right)=
\frac{(2\pi)^{\sfrac12}\mu^{\frac23}g(\mu)}{\phi(\zeta)}\left[
\frac{\Ai(\mu^{\frac43}\zeta)}{\mu^{\frac43}}\,C_\mu(\zeta)+
\Ai'(\mu^{\frac43}\zeta)\,D_\mu(\zeta)\right],\label{a2}
\end{equation}
\begin{equation}
V\left(-\sfrac12\mu^2,\mu t\sqrt{2}\right)=
\frac{2\pi^{\sfrac12}\mu^{\frac13}g(\mu)\phi(\zeta)}{\Gamma(\sfrac12+\sfrac12\mu^2)}\left[
\Bi(\mu^{\frac43}\zeta)\,A_\mu(\zeta)+
\frac{\Bi'(\mu^{\frac43}\zeta)}{\mu^{\frac83}}\,B_\mu(\zeta)\right],\label{a3}
\end{equation} 
\begin{equation}
V'\left(-\sfrac12\mu^2,\mu t\sqrt{2}\right)=
\frac{(2\pi)^{\sfrac12}\mu^{\frac23}g(\mu)}{\phi(\zeta)\Gamma(\sfrac12+\sfrac12\mu^2)}\left[
\frac{\Bi(\mu^{\frac43}\zeta)}{\mu^{\frac43}}\,C_\mu(\zeta)+
\Bi'(\mu^{\frac43}\zeta)\,D_\mu(\zeta)\right].\label{a4}
\end{equation}
The coefficient functions $A_\mu(\zeta), B_\mu(\zeta), C_\mu(\zeta)$ and $D_\mu(\zeta)$
have the following asymptotic expansions
\begin{equation}
 A_\mu(\zeta)\sim\sum_{s=0}^\infty\,    \frac{a_s(\zeta)}{\mu^{4s}},\quad
   B_\mu(\zeta)\sim\sum_{s=0}^\infty\,    \frac{b_s(\zeta)}{\mu^{4s}},\label{a5}
\end{equation}
\begin{equation}
 C_\mu(\zeta)\sim\sum_{s=0}^\infty\,    \frac{c_s(\zeta)}{\mu^{4s}},\quad
   D_\mu(\zeta)\sim\sum_{s=0}^\infty\,    \frac{d_s(\zeta)}{\mu^{4s}},\label{a6}
\end{equation}
as $\mu\to\infty$, uniformly with respect to $t\ge-1+\delta$, where $\delta$ is a 
small fixed positive number. 
The quantity $\zeta$ is defined by
\begeq
\sfrac23(-\zeta)^{\sfrac32}&=\eta(t),\quad-1<t\le1,\quad (\zeta\le0)\\
\sfrac23\zeta^{\sfrac32}&= \xi(t),\quad1\le t,\quad (\zeta\ge0)\label{a7}
\endeq
where $\eta,\xi$ follow from (2.26), (2.2), respectively, and
\begin{equation}
\phi(\zeta)=\left(\frac{\zeta}{t^2-1}\right)^{\sfrac14}.\label{a8}
\end{equation}
The function $\zeta(t)$ is real for $t>-1$ and analytic at $t=1$. We can invert  
$\zeta(t)$ into $t(\zeta)$, and obtain
$$t=1+2^{-1/3}\zeta-\sfrac1{10}2^{-2/3}\zeta^2+\sfrac{11}{700}\zeta^3+\ldots\ .$$
The function $g(\mu)$ has the expansion given in (2.6) and the
coefficients $a_s(\zeta), b_s(\zeta)$ are given by
\begin{equation}
a_s(\zeta)=\sum_{m=0}^{2s}\beta_m\zeta^{-\sfrac32m}{\cal A}_{2s-m}(t)\quad
\sqrt{{\zeta}}b_s(\zeta)=-\sum_{m=0}^{2s+1}\alpha_m\zeta^{-\sfrac32m}{\cal
A}_{2s-m+1}(t),\label{a9}
\end{equation} 
where ${\cal A}_s(t)$ are used in (2.1), $\alpha_0=1$ and
\begin{equation}
\alpha_m=\frac{(2m+1)\,(2m+3)\cdots(6m-1)}{m!\,(144)^m},\quad
\beta_m=-\frac{6m+1}{6m-1}\,a_m.\label{a10}
\end{equation}
A recursion for $\alpha_m$ reads:
$$\alpha_{m+1}=\alpha_m\,\frac{(6m+5)(6m+3)(6m+1)}{144\,(m+1)(2m+1)},\quad
m=0,1,2,\ldots\ .$$
The numbers $\alpha_m, \beta_m$ occur in the asymptotic expansions of the Airy
functions, and the relations in
\parref{a9} follow from solving \parref{a1} and \parref{a3} for $A_\mu(\zeta)$ and
$B_\mu(\zeta)$, expanding the Airy functions (assuming that
$\zeta$ is bounded away from 0) and by using (2.1) and a similar result for $V(a,z)$ 
(O-(11.16) and O-(2.12)).

For negative values of $\zeta$ (i.e., $-1<t<1$) we can use (O-(13.4))
\begeq
a_s(\zeta)&=(-1)^s\sum_{m=0}^{2s}\beta_m(-\zeta)^{-\sfrac32m}
\widetilde{{\cal A}}_{2s-m}(t),\\
\sqrt{{-\zeta}}b_s(\zeta)&=(-1)^{s-1}\sum_{m=0}^{2s+1}\alpha_m(-\zeta)^{-\sfrac32m}
\widetilde{{\cal A}}_{2s-m+1}(t),\label{a11}
\endeq
where
$$\widetilde{{\cal A}}_s(t)=\frac{u_s(t)}{(1-t^2)^{\frac32s}}.$$

The functions $C_\mu(\zeta)$ and $D_\mu(\zeta)$ of \parref{a2} and \parref{a4} are given by
\begin{equation}
C_\mu(\zeta)=\chi(\zeta)A_\mu(\zeta)+A_\mu'(\zeta)+\zeta B_\mu(\zeta),\quad
D_\mu(\zeta)=A_\mu(\zeta)+\frac1{\mu^4}
\left[\chi(\zeta)B_\mu(\zeta)+B_\mu'(\zeta)\right].
\label{a12}
\end{equation}
The coefficients $c_s(\zeta)$ and $d_s(\zeta)$ in \parref{a6} are given by
\begin{equation}
c_s(\zeta)=\chi(\zeta)a_s(\zeta)+a'_s(\zeta)+\zeta b_s(\zeta),\quad
d_s(\zeta)=a_s(\zeta)+\chi(\zeta)b_{s-1}(\zeta)+b'_{s-1}(\zeta),\label{a13}
\end{equation}
where
\begin{equation}
\chi(\zeta)=\frac{\phi'(\zeta)}{\phi(\zeta)}=\frac{1-2t[\phi(\zeta)]^6}{4\zeta},\label{a14}
\end{equation}
with $\phi(\zeta)$ given in \parref{a8}. Explicitly,
\begin{equation}
\frac1{\sqrt{{\zeta}}}c_s(\zeta)=-\sum_{m=0}^{2s+1}\beta_m\zeta^{-\sfrac32m}{\cal B}_{2s-m+1}(\tau)\quad
d_s(\zeta)=-\sum_{m=0}^{2s}\alpha_m\zeta^{-\sfrac32m}{\cal B}_{2s-m}(\tau),\label{a15}
\end{equation}
where ${\cal B}_s(\tau)=v_s(t)/(t^2-1)^{\sfrac32s}$, with $v_s(t)$ 
defined in (2.25). Other versions of \parref{a15} are needed 
for negative values of $\zeta$,i.e., if $-1<t<1$; see \parref{a11}.
\subsection{The case $\quad a\ll0, z\le0$.}%
Near the other turning point $t=-1$ we can use the representations (O-(9.7))
\begeq
&U\left(-\sfrac12\mu^2,-\mu t\sqrt{2}\right)=
2\pi^{\sfrac12}\mu^{\frac13}g(\mu)\phi(\zeta)\Biggl[\sin\left(\sfrac12\pi\mu^2\right)
\Biggl\{\Ai(\mu^{\frac43}\zeta)A_\mu(\zeta)+\\
&\quad\Frac{\Ai'(\mu^{\frac43}\zeta)}{\mu^{\frac83}}B_\mu(\zeta)\Biggr\}
+\cos\left(\sfrac12\pi\mu^2\right)\Biggl\{\Bi(\mu^{\frac43}\zeta)
A_\mu(\zeta)+
\Frac{\Bi'(\mu^{\frac43}\zeta)}{\mu^{\frac83}}
B_\mu(\zeta)\Biggr\}\Biggr],\label{a16}
\endeq 
as $\mu\to\infty$, uniformly with respect to $t\ge-1+\delta$,
where $\delta$ is a small fixed positive number. 
Expansions for $V(a,z)$ follow from \parref{a1} and \parref{a16} and the
second relation in (1.7). Results for the derivatives
of $U(a,z)$ and $V(a,z)$ follow easily from the earlier results.

\subsection{Modified forms of Olver's Airy-type expansions}%
Modified versions of the Airy-type expansions \parref{a1} -- \parref{a4} can also be
given. In the case of the expansions in terms of elementary functions
our main motivation for introducing modified expansions was the double
asymptotic property of these expansions. In the case of the Airy-type
expansions the interesting domains for the parameter $t$, from a numerical
point of view, are finite domains that contain the turning points $\pm1$.
So, considering the expansions given so far, there is no need to have 
Airy-type expansions  with the double asymptotic property; if
$\mu$ remains finite and $|t|\gg1$ we can use the expansions in terms of
elementary functions.  However, we have another interest in modified
expansions in the case of Airy-type expansions. We explain this by
first discussing a few properties of the coefficient functions 
$A_\mu(\zeta), B_\mu(\zeta), C_\mu(\zeta)$ and $D_\mu(\zeta)$.

By using the Wronskian relation (1.6) we can verify the relation
\begin{equation}
A_\mu(\zeta)D_\mu(\zeta)-\frac1{\mu^4}B_\mu(\zeta)C_\mu(\zeta)=
\frac{\Gamma(\sfrac12+\sfrac12\mu^2)}{2\mu\sqrt{{\pi}}g^2(\mu)},\label{a17}
\end{equation}
where $g(\mu)$ is defined by means of an asymptotic expansion given in (2.6).
By using the differential equation (O-(7.2))
\begin{equation}
\frac{d^2W}{d\zeta^2}=\left[\mu^4\zeta+\Psi(\zeta)\right]W,\label{a18}
\end{equation}
where
\begin{equation}
\Psi(\zeta)=\frac5{16\zeta^2}-\frac{(3t^2+2)\zeta}{4(t^2-1)^3}=
2^{1/3}\left[-\sfrac9{280}+\sfrac7{150}2^{-1/3}\zeta
-\sfrac{1359}{26950}2^{-2/3}\zeta^2+\sfrac{196}{8125}\zeta^3\ldots\right],\label{a19}
\end{equation}
we can derive the following system of equations for the functions $A_\mu(\zeta),
B_\mu(\zeta)$:
\begeq
A''+2\zeta B'+B-\Psi(\zeta)A&=0,\\
B''+2\mu^4A'-\Psi(\zeta)B&=0,\label{a20}
\endeq
where primes denote differentiation with respect to $\zeta$. 
A Wronskian for this system follows by eliminating the terms with $\Psi(\zeta)$.
This gives
$$2\mu^4A'A+AB''-A''B-2\zeta B'B-B^2=0,$$
which can be integrated:
\begin{equation}
\mu^4A_\mu^2(\zeta)+A_\mu(\zeta)B_\mu'(\zeta)-A_\mu'(\zeta)B_\mu(\zeta)-\zeta B_\mu^2(\zeta)=
\mu^4\frac{\Gamma(\sfrac12+\sfrac12\mu^2)}{2\mu\sqrt{{\pi}}g^2(\mu)},\label{a21}
\end{equation}
where the quantity on the right-hand side follows from \parref{a17} and \parref{a12}.
It has the expansion 
\begin{equation}
\mu^4\left[1-\frac1{576\mu^4}+\frac{2021}{2488320\mu^8}+\ldots\right],\label{a22}
\end{equation}
as follows from O-(2.22) and O-(5.21).

As mentioned before, the interesting domain of the Airy-type expansions given
in this section is the domain that contains the turning point $t=1$, or
$\zeta=0$.  The representations of the coefficients of the expansions given in
\parref{a9} cannot be used in numerical algorithms when $|\zeta|$ is small, unless
we expand all relevant coefficients in powers of $\zeta$. This is one way how to handle
this  problem numerically; see \cite{temna}. In that paper we have 
discussed another method that is based on 
a system like \parref{a20}, with applications to Bessel functions.
In that method the functions $A_\mu(\zeta)$ and $B_\mu(\zeta)$ are expanded 
in powers of $\zeta$, for sufficiently small values of $|\zeta|$, 
say $|\zeta|\le1$, and the Maclaurin coefficients are 
computed from \parref{a20} by recursion. A normalizing relation (the analogue of
\parref{a21}) plays a crucial role in that algorithm. The method works quite well
for relatively small values of a parameter (the order of the Bessel functions)
that is the analogue of $\mu$. 

When we want to use this algorithm for the present case only large values of
$\mu$ are allowed because the function $g(\mu)$ that is used in
\parref{a1}--\parref{a4} and \parref{a21} is only defined for large values of $\mu$. 
For this reason we give the modified versions of Olver's  Airy-type expansions.
The modified versions are more
complicated  than the Olver's expansions, because the analogues
of the series in
\parref{a5} and \parref{a6} are in powers of
$\mu^{-2}$,  and not in powers of $\mu^{-4}$. Hence, when we use these
series for numerical computations we need more coefficients in the modified
expansions, which is certainly not desirable from a numerical point of view,
given the complexity of the coefficients in Airy-type expansions. However, in
the algorithm based on Maclaurin expansions of the analogues of 
the coefficient functions  $A_\mu(\zeta), B_\mu(\zeta), C_\mu(\zeta)$ and $D_\mu(\zeta)$
this point is of minor concern.

The modified expansions are the following:

\begin{equation}
U\left(-\sfrac12\mu^2,\mu t\sqrt{2}\right)=
\frac{\Gamma(\sfrac12+\sfrac12\mu^2)\phi(\zeta)}{\mu^{\frac23}h(\mu)}\left[
\Ai(\mu^{\frac43}\zeta)\,F_\mu(\zeta)+
\frac{\Ai'(\mu^{\frac43}\zeta)}{\mu^{\frac83}}\,G_\mu(\zeta)\right],\label{a23}
\end{equation} 
\begin{equation}
V\left(-\sfrac12\mu^2,\mu t\sqrt{2}\right)=
\frac{\phi(\zeta)}{\mu^{\frac23}h(\mu)}\left[
\Bi(\mu^{\frac43}\zeta)\,F_\mu(\zeta)+
\frac{\Bi'(\mu^{\frac43}\zeta)}{\mu^{\frac83}}\,G_\mu(\zeta)\right].\label{a24}
\end{equation} 
The functions $F_\mu(\zeta)$ and $G_\mu(\zeta)$ have the
following asymptotic expansions
\begin{equation}
 F_\mu(\zeta)\sim\sum_{s=0}^\infty\,    \frac{f_s(\zeta)}{\mu^{2s}},\quad
   G_\mu(\zeta)\sim\sum_{s=0}^\infty\,    \frac{g_s(\zeta)}{\mu^{2s}}.\label{a25}
\end{equation}
The quantity $\zeta$ and the functions $\phi(\zeta)$ and 
$h(\mu)$ are as in Section 3.1.
Comparing \parref{a23}, \parref{a24} with \parref{a1}, \parref{a3} we conclude that
\begin{equation}
F_\mu(\zeta)=H(\mu)A_\mu(\zeta),\quad
G_\mu(\zeta)=H(\mu)B_\mu(\zeta),\quad H(\mu)=\frac{2\sqrt{{\pi}}\mu
g(\mu)h(\mu)}{\Gamma(\sfrac12+\sfrac12\mu^2)}.\label{a26}
\end{equation}
The function $H(\mu)$ can be expanded (see O-(2.22), O-(2.27), O-(6.2) and 
(2.6))
\begin{equation}
H(\mu)\sim 1+\sfrac12\sum_{s=1}^\infty(-1)^s
\frac{\gamma_s}{(\sfrac12\mu^2)^s},\label{a27}
\end{equation}
where $\gamma_s$ are the coefficients in the gamma function expansions
\begin{equation}
\Gamma(\sfrac12+z)\sim\sqrt{{2\pi}}e^{-z}z^z\sum_{s=0}^\infty\frac{\gamma_s}{z^s},\quad
\frac1{\Gamma(\sfrac12+z)}\sim\frac{e^{z}z^{-z}}{\sqrt{{2\pi}}}
\sum_{s=0}^\infty(-1)^s\frac{\gamma_s}{z^s}.\label{a28}
\end{equation}
The first few coefficients are
$$\gamma_0=1,\quad \gamma_1=-\frac1{24},\quad \gamma_2=\frac1{1152},\quad
\gamma_3=\frac{1003}{414720}.$$

The second expansion in \parref{a28} can be used in \parref{a26} to find relations
between the  coefficients $a_s(\zeta)$ and $b_s(\zeta)$ of \parref{a5} and of 
$f_s(\zeta)$ and $g_s(\zeta)$ of \parref{a25}. That is,
$$f_0(\zeta)=1,\quad f_1(\zeta)=\sfrac1{24},\quad f_2(\zeta)=a_1(\zeta)+\sfrac1{576},
\quad f_3(\zeta)=\sfrac1{24}a_1(\zeta)-\sfrac{1003}{103680},$$
$$g_0(\zeta)=b_0(\zeta),\  g_1(\zeta)=\sfrac1{24}b_0(\zeta),\ 
g_2(\zeta)=b_1(\zeta)+\sfrac1{576}b_0(\zeta),\ 
g_3(\zeta)=\sfrac1{24}b_1(\zeta)-\sfrac{1003}{103680}b_0(\zeta).$$ 
The coefficients $f_s(\zeta), g_s(\zeta)$
can also be expressed in terms of the coefficients $\phi_s(\tau)$ that are introduced in
(2.9) by deriving the analogues of \parref{a9}.

 The system of equations \parref{a20} remains the same:
\begeq
F''+2\zeta G'+G-\Psi(\zeta)F&=0,\\
G''+2\mu^4F'-\Psi(\zeta)G&=0,\label{a29}
\endeq
and the Wronskian relation becomes
\begin{equation}
\mu^4F_\mu^2(\zeta)+F_\mu(\zeta)G_\mu'(\zeta)-F_\mu'(\zeta)G_\mu(\zeta)-\zeta G_\mu^2(\zeta)=
\mu^4\frac{2\sqrt{{\pi}}\mu h^2(\mu)}{\Gamma(\sfrac12+\sfrac12\mu^2)}.\label{a30}
\end{equation}
The right-hand side has the expansion (see \parref{a28} and (2.7))
$\mu^4\sum_{s=0}^\infty(-1)^s \gamma_s/(\sfrac12\mu^2)^s$.
Observe that (\parref{a30} is an exact relation, whereas \parref{a21} contains the function
$g(\mu)$, of which only an asymptotic expansion is available.
\subsection{Numerical aspects of the Airy-type expansions}%
In \cite{temna}, (Section 4), we solved the system \parref{a29} (for the case of Bessel
functions) by substituting Maclaurin series of $F(\zeta), G(\zeta)$ and $\Psi(\zeta)$. That is,
we wrote
$$F(\zeta)=\sum_{n=0}^\infty c_n(\mu)\zeta^n,\quad G(\zeta)=\sum_{n=0}^\infty d_n(\mu)\zeta^n,\quad \Psi(\zeta)=\sum_{n=0}^\infty \psi_n\zeta^n,$$
where the coefficients $\psi_n$ can be considered as known (see \parref{a19}),
and we substituted the expansions in
\parref{a29}. This gives for $n=0, 1, 2,\ldots$ the recursion relations
\begeq
(n+2)(n+1)c_{n+2}+(2n+1)d_n=\rho_n,\quad &\rho_n={\displaystyle\sum_{k=0}^n} \psi_kc_{n-k},
\morespace\\
(n+2)(n+1)d_{n+2}+2\mu^4(n+1)c_{n+1}=\sigma_n,\quad &\sigma_n={\displaystyle\sum_{k=0}^n}
\psi_kd_{n-k}.\label{a31}
\endeq
If $\mu$ is large, the recursion relations cannot be solved in forward direction,
because of numerical instabilities. For the Bessel function case we have shown that we
can solve the system by iteration and backward recursion. The relation in \parref{a30}
can be used for normalization of the coefficients in the backward recursion scheme.

For details we refer to 
\cite{temna}. The present case is identical to the case of the Bessel functions;
only the function $\Psi(\zeta)$ is different, and instead of $\mu^2$ in \parref{a31} we had
the order $\nu$ of the Bessel functions.
\section{Expansions from integral representations}%
The expansions developed by Olver, of which some are given in the previous sections,  
are all valid if
$|a|$ is large. For several cases we gave modified expansions that hold if at least
one of the two parameters $a, z$ is large and we have indicated the relations
between Olver's expansions and the new
expansions. The modified expansions have in fact a double
asymptotic property. Initially, we derived these expansions by using integral
representations of the parabolic cylinder functions, and later we found the
relations with Olver's expansions.  In this section we explain how some of the
modified expansions can be obtained from the integrals that define
$U(a,z)$ and $V(a,z)$. Again we only consider real values of the parameters.
\subsection{Expansions in terms of elementary functions by using integrals}%
\subsubsection{The case $\quad a\ge0, z\ge0;\quad a+z\gg0$.}%
We start with the well-known integral representation
\begin{equation}
U(a,z)=\frac{e^{-\sfrac14z^2}}{\Gamma(a+\sfrac12)}
\intp w^{a-\sfrac12} e^{-\sfrac12w^2-zw}\,dw,\quad a>-\sfrac12\label{i1}
\end{equation}
which we write in the form
\begin{equation}
U(a,z)=\frac{z^{a+\sfrac12}\,e^{-\sfrac14z^2}}{\Gamma(a+\sfrac12)}\intp w^{-\sfrac12}
e^{-z^2\phi(w)}\, dw,\label{i2}
\end{equation}
where
\begin{equation}
\phi(w)=w+\sfrac12w^2-\lambda\ln w,\quad \lambda=\frac{a}{z^2}.\label{i3}
\end{equation}
The positive saddle point $w_0$ of the integrand in \parref{i3} is computed
from
\begin{equation}
\frac{d\phi(w)}{dw}=\frac{w^2+w-\lambda}{w}=0,\label{i4}
\end{equation}
giving
\begin{equation}
 w_0=\sfrac12\left[\sqrt{{1+4\lambda}}-1\right].\label{i5}
\end{equation}
We consider $z$ as the large parameter. When $\lambda$ is bounded away from 0
we can use Laplace's method (see \cite{olver} or \cite{wong}). 
When $a$ and $z$ are such that
$\lambda\to0$ Laplace's method cannot be applied. However, we can use a
method given in \cite{temlap} that allows small values
of $\lambda$. 

To obtain a standard form for this Laplace-type integral, we
transform $w\to t$ (see \cite{temq})
by writing
\begin{equation}
\phi(w)=t-\lambda\ln t+A,\label{i6}
\end{equation}
where $A$ does not depend on $t$ or $w$, and we prescribe that $w=0$
should correspond with $t=0$ and $w=w_0$  with $t=\lambda$, the saddle point
in the $t-$plane.

This gives
\begin{equation}
U(a,z)=\frac{z^{a+\sfrac12}\,e^{-\sfrac14z^2-Az^2}}{(1+4\lambda)^{\frac14}\Gamma(a+\sfrac12)}\intp
t^{a-\sfrac12}e^{-z^2t} f(t)\,dt,\label{i7}
\end{equation}
where
\begin{equation}
f(t)=(1+4\lambda)^{\sfrac14}\sqrt{{\frac tw }}\,\frac{dw}{dt}=(1+4\lambda)^{\sfrac14}
\sqrt{{\frac wt }}\,\frac{t-\lambda}{w^2+w-\lambda}.\label{i8}
\end{equation}
By normalizing with the quantity $(1+4\lambda)^{\sfrac14}$ we obtain $f(\lambda)=1$,
as can be verified from \parref{i8} and a limiting process (using
l'H\^opital's rule). The quantity $A$ is given by
\begin{equation}
A=\sfrac12w_0^2+w_0-\lambda\ln w_0-\lambda+\lambda\ln\lambda.\label{i9}
\end{equation}

A first uniform expansion can be obtained by writing
\begin{equation}
f(t)=\sum_{n=0}^\infty a_n(\lambda)(t-\lambda)^n.\label{i10}
\end{equation}
Details on the computation of $a_n(\lambda)$ will be given in the appendix.

By substituting \parref{i10} into \parref{i7} we obtain
\begin{equation}
U(a,z)\sim \frac{e^{-\sfrac14z^2-Az^2}}{z^{a+\sfrac12}\,(1+4\lambda)^{\frac14}}
\sum_{n=0}^\infty a_n(\lambda) P_n(a) z^{-2n},\label{i11}
\end{equation}
where
\begin{equation}
P_n(a)=\frac{z^{2a+2n+1}}{\Gamma(a+\sfrac12)}\intp
t^{a-\sfrac12}\,e^{-z^2t}(t-\lambda)^n\,dt,\quad n=0,1,2,\ldots\ .\label{i12}
\end{equation}
The $P_n(a)$ are polynomials in $a$. They follow the recursion
relation
$$P_{n+1}(a)=(n+\sfrac12)P_n(a)+an P_{n-1}(a), \quad
n=0,1,2,\ldots,$$ 
with initial values
$$ P_0(a)=1,\quad P_1(a)=\sfrac12.$$

We can obtain a second expansion
\begin{equation}
U(a,z)\sim 
\frac{e^{-\sfrac14z^2-Az^2}}{z^{a+\sfrac12}\,(1+4\lambda)^{\frac14}}
\sum_{k=0}^\infty \frac{f_k(\lambda)}{z^{2k}},\label{i13}
\end{equation}
with the property that in the series the parameters $\lambda$ and $z$ are
separated, by introducing a sequence of functions $\{f_k\}$ with
$f_0(t)=f(t)$ and by defining
\begin{equation}
f_{k+1}(t)=\sqrt{{t}}\,\frac{d\ }{dt}\left[\sqrt{{t}}\frac{f_k(t)-f_k(\lambda)}{t-\lambda}\right],\quad
k=0,1,2,\ldots\ .\label{i14}
\end{equation}
The coefficients $f_k(\lambda)$ can be expressed in terms of the
coefficients $a_n(\lambda)$ defined in \parref{i10}. To verify this, we write
\begin{equation}
f_k(t)=\sum_{n=0}^\infty a_n^{(k)}(\lambda)(t-\lambda)^n,\label{i15}
\end{equation}
and by substituting this in \parref{i14} it follows that
\begin{equation}
a_n^{(k+1)}(\lambda)=\lambda(n+1) a_{n+2}^{(k)}(\lambda)+\left(n+\sfrac12\right)
a_{n+1}^{(k)}(\lambda),\quad k\ge0,\quad
n\ge0.\label{i16}
\end{equation}
Hence, the coefficients $f_k(\lambda)$ of \parref{i13} are given by
\begin{equation}
f_k(\lambda)=a_0^{(k)}(\lambda),\quad k\ge0.\label{i17}
\end{equation}
We have
\begeq
f_0(\lambda)&= 1,\\
f_1(\lambda)&=  \sfrac12\left[a_1(\lambda) + 2 \lambda a_2(\lambda)\right],\\
f_2(\lambda)&=  \sfrac14\left[12 \lambda^2 a_4(\lambda)
+ 14 \lambda a_3(\lambda) + 3 a_2(\lambda)\right],\\
f_3(\lambda)&=  \sfrac18\left[120 \lambda^3 a_6(\lambda) + 220 \lambda^2 a_5(\lambda) +
116 \lambda a_4(\lambda) + 15 a_3(\lambda)\right].\label{i18}
\endeq
Explicitly,
\begeq
f_0(\lambda)&=   1,\morespace\\
f_1(\lambda)&=   -\Frac\rho{24}  (20\sigma^2-10\sigma-1),\morespace\\
f_2(\lambda)&=   \Frac{\rho^2}{1152}(6160\sigma^4-6160\sigma^3+924\sigma^2+20\sigma+1),\morespace\\
f_3(\lambda)&=-\Frac{\rho^3}{414720}(27227200\sigma^6-40840800\sigma^5+16336320\sigma^4-\morespace\\
        &   \quad\quad 1315160\sigma^3-8112\sigma^2+2874\sigma+1003),\label{i19}
\endeq
where
\begin{equation}
\sigma=\sfrac12\left[1+\frac{z}{\sqrt{{4a+z^2}}}\right],\quad
\rho=\frac{(2\sigma-1)^2}{\sigma}=\frac{2z^2}{\sqrt{{4a+z^2}}(z+\sqrt{{4a+z^2}})}.\label{i20}
\end{equation}
We observe that $f_k(\lambda)$ is a polynomial of degree $2k$ in $\sigma$
multiplied with $\rho^k$.

If $a$ and $z$ are positive then $\sigma\in[0,1]$. Furthermore, the sequence
$\{\rho^k/z^{2k}\}$ is an asymptotic scale when one or both
parameters $a$ and $z$ are large.
The expansion in \parref{i13} is valid for $z\to\infty$ and holds uniformly
for $a\ge0$. It has a double asymptotic property in the sense that it is
also valid as $a\to\infty$, uniformly with respect to $z\ge 0$. As follows
from the coefficients given in \parref{i19} and relations to be
given later, we can indeed let $z\to0$ in the expansion.
%
%
%

The expansion in \parref{i13} can be obtained by using an integration by
parts procedure. We give a few steps in this method.
Consider the integral
\begin{equation}
F_a(z)=\frac1{\Gamma(a+\sfrac12)}
\intp t^{a-\sfrac12}e^{-z^2t} f(t)\,dt,\label{i21}
\end{equation}
We have (with $\lambda=a/z^2$)
\begin{eqnarray*}
F_a(z)
&=\Frac{f(\lambda)}{\Gamma(a+\sfrac12)}
\intp t^{a-\sfrac12}e^{-z^2t} \,dt +
\Frac{1}{\Gamma(a+\sfrac12)}
\intp t^{a-\sfrac12}e^{-z^2t}[f(t)-f(\lambda)] \,dt\\
&= z^{-2a-1}f(\lambda)  -
\Frac{1}{z^2\Gamma(a+\sfrac12)}
\intp t^{\sfrac12}\Frac{[f(t)-f(\lambda)]}{t-\lambda} \,de^{-z^2(t-\lambda\ln t)}\\
&= z^{-2a-1}f(\lambda) +  \Frac1{z^2\Gamma(a+\sfrac12)}
\intp t^{a-\sfrac12}e^{-z^2t} f_1(t)\,dt,
\end{eqnarray*}
where $f_1$ is given in \parref{i14} with $f_0=f$. Repeating this
procedure we obtain \parref{i13}. More details on this method and
proofs of the asymptotic nature of the expansions \parref{i11} and
\parref{i13} can be found in our earlier papers. We concentrate on
expansion \parref{i13} because \parref{i11} cannot be compared
with Olver's expansions.

To compare \parref{i13} with Olver's expansion (2.16), we write
\begin{equation}
a=\sfrac12\mu^2,\quad z=\mu\sqrt{2} t.\label{i22}
\end{equation}
Then the parameters $\sigma$ and $\rho$ defined in \parref{i20} become
\begin{equation}
\sigma=\sfrac12\left[1+\frac{t}{\sqrt{{1+t^2}}}\right]=\widetilde\tau+1,\quad
\rho=\frac{2t^2}{\sqrt{{1+t^2}}(t+\sqrt{{1+t^2}})},\label{i23}
\end{equation}
where $\widetilde\tau$ is given in (2.32). 
After a few manipulations we write \parref{i13} in the form (cf.(2.29)) 
\begin{equation}
U\left(\sfrac12\mu^2,\mu t\sqrt{2}\right)=
\frac{\widetilde h(\mu)\,e^{-\mu^2{\widetilde\xi}}}
{(t^2+1)^{\sfrac14}}\,\widetilde F_\mu(z),\quad
\widetilde F_\mu(z)\sim\sum_{k=0}^\infty(-1)^k\frac{\widetilde\phi_k(\sigma)}{\mu^{2k}},\label{i24}
\end{equation}
where 
\begin{equation}
\widetilde\xi=\sfrac12\left[t\sqrt{{1+t^2}}+\ln\left(t+\sqrt{{1+t^2}}\right)\right],
\label{i25}
\end{equation}
\begin{equation}
h(\mu)=e^{\sfrac14\mu^2}\,\mu^{-\sfrac12\mu^2-\sfrac12}\,
2^{\sfrac14\mu^2-\sfrac14}\label{i26}
\end{equation}
and
\begin{equation}
\widetilde\phi_k(\sigma)=\frac{(-1)^k}{(2t^2)^{k}}f_k(\lambda).\label{i27}
\end{equation}
Explicitly,
\begeq
\widetilde\phi_0(\sigma)&=   1,\morespace\\
\widetilde\phi_1(\sigma)&=   \Frac{1-\sigma}{12}  (20\sigma^2-10\sigma-1),\morespace\\
\widetilde\phi_2(\sigma)&=   
\Frac{(1-\sigma)^2}{288}(6160\sigma^4-6160\sigma^3+924\sigma^2+20\sigma+1),\morespace\\
\widetilde\phi_3(\sigma)&=  
\Frac{(1-\sigma)^3}{51840}(27227200\sigma^6-40840800\sigma^5+16336320\sigma^4-\morespace\\
       &   \quad\quad 1315160\sigma^3-8112\sigma^2+2874\sigma+1003),\label{i28}
\endeq
where $\sigma$ is given in \parref{i23}.
Comparing \parref{i24} with (2.29) we obtain 
$\widetilde\phi_k(\sigma)=\phi_k(\widetilde\tau), k\ge0$, because $\sigma=1+\widetilde\tau$.

\subsubsection{The case $\quad a\ge0, z\le0;\quad a-z\gg0$.}
To derive the first expansion in (2.34) we use the contour integral
\begin{equation}
U(a,-z)=\frac{\sqrt{{2\pi}}e^{\sfrac14z^2}}{\Gamma(a+\sfrac12)}H_a(z),\quad
H_a(z)=\frac{\Gamma(a+\sfrac12)}{2\pi
i}\int_{{\cal C}} e^{zs+\sfrac12 s^2}s^{-a-\sfrac12}\,ds,\label{i40}
\end{equation}
where ${{\cal C}}$ is a vertical line in the half plane $\Re s>0$. This integral
can be transformed into a standard form that involves the same mapping
as in the previous subsection. We first write (by transforming via $s=zw$)
\begin{equation}
H_a(z)=\frac{z^{\sfrac12-a}\Gamma(a+\sfrac12)}{2\pi i}\int_{{\cal C}}
e^{z^2(w+\sfrac12 w^2)}w^{-a-\sfrac12}\,dw
=\frac{z^{\sfrac12-a}\Gamma(a+\sfrac12)}{2\pi i}\int_{{\cal C}}
e^{z^2\phi(w)}\,\frac{dw}{\sqrt{{w}}},\label{i41}
\end{equation}
where $\phi(w)$ is defined in \parref{i3}.
By using the transformation given in \parref{i6} it follows that
\begin{equation}
H_a(z)=\frac{z^{\sfrac12-a}\Gamma(a+\sfrac12)e^{Az^2}}{2\pi i}\int_{{\cal C}}
e^{z^2t}t^{-a-\sfrac12}\,f(t)\,dt.\label{i42}
\end{equation}
The integration by parts method used for \parref{i21} gives the
expansion (see \cite{temlap})
\begin{equation}
H_a(z)\sim\frac{z^{a} e^{Az^2}}
{(4a+z^2)^{\frac14}}\sum_{k=0}^\infty(-1)^k\frac{f_k(\lambda)}{z^{2k}},\label{i43}
\end{equation}
where the $f_k(\lambda)$ are the same as in \parref{i13}; see also \parref{i18}.
This gives the first expansion of (2.34).

\noindent
{\bf Remark 4.1.\ }
The first result in (2.34) can also be obtained by using \parref{i1}
with $z<0$. The integral for $U(a,-z)$ can be written as in
\parref{i2}, now with $\phi(w)=\sfrac12w^2-w-\ln\lambda,\ \lambda = a/z^2$. 
In this case the relevant saddle point at $w_0=(1+\sqrt{{1+4\lambda}})/2$ is 
always inside the interval
$[1,\infty)$ and the standard method of Laplace can be used. The same
expansion will be obtained with the same structure and coefficients as
in (2.34), because of the unicity of Poincar\'e-type asymptotic
expansions. See also Section 4.1.4 where Laplace's method will be used for an
integral that defines $V(a,z)$.
\subsubsection{The case $\quad a\le0, z>2\sqrt{{-a}},\quad-a+z\gg0$.}%
Olver's starting point (2.1) can also be obtained from an integral. 
Observe that \parref{i1} is not valid for $a\le-\sfrac12$. We take as
integral (see \cite{abst}, p. 687, 19.5.1)
\begin{equation}
U(-a,z)=\frac{\Gamma(\sfrac12+a)}{2\pi i}e^{-\sfrac14z^2}\int_\alpha
e^{zs-\sfrac12s^2}\,s^{-a-\sfrac12}\,ds,\label{i46}
\end{equation}
where $\alpha$ is a contour that encircles the negative $s-$axis in positive
direction. Using a transformation we can write this in the form (cf.
\parref{i2})
\begin{equation}
U(-a,z)=\frac{\Gamma(\sfrac12+a)}{2\pi i}z^{\sfrac12-a}e^{-\sfrac14z^2}\int_\alpha
e^{\phi(w)}\,w^{-\sfrac12}\,dw,\label{i47}
\end{equation}
where 
\begin{equation}
\phi(w)=w-\sfrac12w^2-\lambda\ln w,\quad \lambda=\frac{a}{z^2}.\label{i48}
\end{equation}
The relevant saddle point is now given by
\begin{equation}
w_0=\sfrac12\left(1-\sqrt{{1-4\lambda}}\right),\quad 0<\lambda<\sfrac14.\label{i49}
\end{equation}

When $\lambda\to0$ the standard saddle point method is not applicable, and
we can again use the methods of
our earlier papers (\cite{temlap} and \cite{temq}) and transform
\begin{equation}
\phi(w)=t-\lambda\ln t+A,\label{i50}
\end{equation}
where the points at $-\infty$ in the $w-$ and $t-$ plane should correspond, 
and $w=w_0$ with $t=\lambda$.  We obtain
\begin{equation}
U(-a,z)=\frac{\Gamma(\sfrac12+a)}{(1-4\lambda)^{\frac14}\,2\pi i}
z^{\sfrac12-a}e^{-\sfrac14z^2+z^2A}\int_\alpha
e^{z^2 t}\,t^{-a-\sfrac12}\,f(t)\,dt,\label{i51}
\end{equation}
where $\alpha$ is a contour that encircles the negative $t-$axis in positive
direction and
\begin{equation}
f(t)=(1-4\lambda)^{\sfrac14}\sqrt{{\frac tw}}\,\frac{dw}{dt}=(1-4\lambda)^{\sfrac14}\,
\sqrt{{\frac wt}}\,\frac{t-\lambda}{w-w^2-\lambda}.\label{i52}
\end{equation}
Expanding $f(t)$ as in \parref{i10}, and computing $f_k(\lambda)$ as in the
procedure that yields the relations in \parref{i18}, we find that the same values 
$f_k(\lambda)$ as in \parref{i19}, up to a factor$(-1)^k$ and a different value
of $\tau$ and $\rho$. By using the integration by parts method for 
contour integrals \cite{temlap}, that earlier produced \parref{i43},
we obtain the result
\begin{equation}
U(-a,z)\sim\frac{z^a\,e^{Az^2-\sfrac14z^2}}
{(z^2-4a)^{\frac14}}\,\sum_{k=0}^\infty(-1)^k\frac{f_k(\lambda)}{z^{2k}},\label{i53}
\end{equation}
where the first $f_k(\lambda)$ are given in \parref{i19} with 
\begin{equation}
\sigma=\sfrac12\left[1+\frac{z}{\sqrt{{z^2-4a}}}\right],\quad
\rho=\frac{(2\sigma-1)^2}{\sigma}=\frac{2z^2}{\sqrt{{z^2-4a+}}(z+\sqrt{{z^2-4a}})}.\label{i54}
\end{equation}
This expansion can be written in the form (2.9). 
\subsubsection{The case $\quad a\le0, z<-2\sqrt{{-a}},\quad-a-z\gg0$.}
We use the relation (see (1.7))
\begin{equation}
U(-a,-z)=\sin\pi a\,U(-a,z)+\frac{\pi}{\Gamma(\sfrac12-a)}\,V(-a,z),\label{i58}
\end{equation}
and use the result of $U(-a,z)$ given in \parref{i53} or the form (2.9). 
An expansion for  $V(-a,z)$
in \parref{i58} can be obtained from the integral (see \cite{mil2})
\begin{equation}
V(a,z)=\frac{e^{-\sfrac14z^2}}{2\pi}\int_{\gamma_1\cup\gamma_2}e^{-\sfrac12s^2+zs}s^{a-\sfrac12}
\,ds,\label{i59}
\end{equation}
where $\gamma_1$ and $\gamma_2$ are two horizontal lines, $\gamma_1$ in the upper
half plane $\Im s>0$ and $\gamma_2$ in the lower half plane $\Im s<0$; the 
integration is from $\Re s=-\infty$ to $\Re s=+\infty$. (Observe that when we
integrate on $\gamma_1$ in the other direction (from $\Re s=+\infty$ to $\Re s=\infty$)
the contour ${\gamma_1\cup\gamma_2}$ can be deformed into $\alpha$ of \parref{i46}, and the
integral defines $U(a,z)$, up to a factor.) We can apply Laplace's method
to obtain the expansion  given in (2.14) (see Remark 4.1).
\subsection{The singular points of the mapping \parref{i6}}%
The mapping defined in \parref{i6} is singular at the saddle point
\begin{equation}
w_-=-\sfrac12\left(\sqrt{{1+4\lambda}}+1\right).\label{i60}
\end{equation}
If $\lambda=0$ then $w_-=-1$ and the corresponding $t-$value is $-\sfrac12$.
For large values of $\lambda$ we have the estimate:
\begin{equation}
t(w_-)\sim\lambda\left[-0.2785-\frac{0.4356}{\sqrt{{\lambda}}}\right].
\label{i61}
\end{equation}
This estimate is obtained as follows. The value $t_-=t(w_-)$ is implicitly
defined by equation \parref{i6} with $w=w_-$. This gives
\begeq
t_--\lambda\ln t_--\lambda+\lambda\ln\lambda
&=-\sfrac12\sqrt{{1+4\lambda}}\pm \lambda\pi i+\lambda\ln\Frac{4\lambda}
{(1+\sqrt{{1+4\lambda}})^2}\\
&=\pm\lambda\pi i-2\sqrt{{\lambda}}\left[1+\Frac1{24\lambda}+\bo(\lambda^{-2})\right],\label{i62}
\endeq
as $\lambda\to\infty$. The numerical solution of the  equation
$s-\ln s-1=\pm\pi i$
is given by
$s_\pm=0.2785\cdots e^{\mp\pi i}$.
This gives the leading term in \parref{i62}. The other term follows by a
further simple step.
\subsection{Expansions in terms of Airy functions}%
All results for the modified Airy-type expansions given in Section 3.3 can 
be obtained by using certain loop integrals. 
The integrals in (4.33) and (4.43)
can be used for obtaining (3.23)
and (3.24), respectively.
The method is based on replacing $\phi(w)$ in (4.34) by a cubic polynomial,
in order to take into account the influence of both saddle points of $\phi(w)$.
This method is first described in \cite{ches}; 
see also \cite{olver} and \cite{wong}.
\section{Numerical verifications}%
We verify several asymptotic expansions by computing the error in the 
Wronskian relation for the series in the asymptotic expansions. Consider 
Olver's expansions of Section 2.3 for the oscillatory region $-1<t<1$ 
with negative $a$. We verify the relation in (2.28). Denote the left-hand side of
the first line in (2.28) by $W(\mu,t)$. Then we define as the error in the
expansions
\begin{equation}
\Delta(\mu,t):=\left|\Frac{W(\mu,t)}
{1-\frac1{576\mu^4}+\frac{2021}{2488320\mu^8}}-1\right|\label{n1}
\end{equation}
Taking 3 terms in the series of (2.23), (2.24) and (2.27), we obtain
for several values of $\mu$ and $t$ the  results given in Table 5.1.
We clearly see the loss of accuracy when $t$ is close to 1.
Exactly the same results are obtained for negative values of $t$ in this
interval.
$$\vbox{
\begin{tabbing}
 10000\=1000\=10000000\=10000000\=10000000\=10000000\=10000000\kill
 \> $\mu$   \>\ \ \ \ 5 \>\ \ \ 10\ \>\  \ \ 25 \>\ \  \ 50 \> \ \ 100  \\
\ $\ t$ \>  \>          \>          \>          \>          \>           \\
\   .00 \>  \>  .32e-09 \>  .78e-13 \>  .13e-17 \>  .32e-21 \>  .78e-25  \\
\   .10 \>  \>  .26e-09 \>  .63e-13 \>  .11e-17 \>  .26e-21 \>  .63e-25  \\
\   .20 \>  \>  .81e-10 \>  .20e-13 \>  .33e-18 \>  .82e-22 \>  .20e-25  \\
\   .30 \>  \>  .16e-08 \>  .39e-12 \>  .65e-17 \>  .16e-20 \>  .39e-24  \\
\   .40 \>  \>  .88e-08 \>  .22e-11 \>  .36e-16 \>  .89e-20 \>  .22e-23  \\
\   .50 \>  \>  .51e-07 \>  .13e-10 \>  .21e-15 \>  .52e-19 \>  .13e-22  \\
\   .60 \>  \>  .40e-06 \>  .99e-10 \>  .17e-14 \>  .40e-18 \>  .99e-22  \\
\   .70 \>  \>  .53e-05 \>  .13e-08 \>  .22e-13 \>  .54e-17 \>  .13e-20  \\
\   .80 \>  \>  .20e-03 \>  .50e-07 \>  .84e-12 \>  .20e-15 \>  .50e-19  \\
\   .90 \>  \>  .35e-00 \>  .24e-04 \>  .41e-09 \>  .10e-12 \>  .25e-16  \\
\end{tabbing}
}$$

\vspace*{-1.4cm}
\begin{quote}
{\bf Table 5.1.}\quad Relative accuracy 
$\Delta(\mu,t)$ defined in \parref{n1} for the asymptotic series of Section 2.3.
\end{quote}

Next we consider the modified expansions of Section 2.1. Denote the left-hand 
side of   (2.20) by $W(\mu,t)$. Then we define as the error in the
expansions
\begin{equation}
\Delta(\mu,t):=\left|{\sfrac12W(\mu,t)}-1\right|\label{n2}
\end{equation}
When we use the series in (2.9), (2.14), (2.18) and (2.19) with 5 terms, we obtain
the results given in Table 5.2. We observe that the accuracy improves as $\mu$ or $t$
increase. This shows the double asymptotic poperty of the modified 
expansions of Section 2.1.
$$\vbox{
\begin{tabbing}
 10000\=1000\=10000000\=10000000\=10000000\=10000000\=10000000\kill
 \>  $\mu$   \>\ \ \ \ 5 \>\ \ \ 10\ \>\  \ \ 25 \>\ \  \ 50 \> \ \ 100   \\
\ $\ t$ \>  \>          \>         \>         \>         \>          \\
\  \ 1.1   \> \>  .51e-01   \>  .48e-05   \>  .72e-10   \>  .18e-13   \>  .43e-17   \\
\  \ 1.2   \> \>  .39e-04   \>  .79e-08   \>  .13e-12   \>  .32e-16   \>  .78e-20   \\
\  \ 1.3   \> \>  .83e-06   \>  .19e-09   \>  .32e-14   \>  .78e-18   \>  .19e-21   \\
\  \ 1.4   \> \>  .56e-07   \>  .13e-10   \>  .23e-15   \>  .55e-19   \>  .13e-22   \\
\  \ 1.5   \> \>  .71e-08   \>  .17e-11   \>  .29e-16   \>  .70e-20   \>  .17e-23   \\
\  \ 2.0   \> \>  .10e-10   \>  .25e-14   \>  .43e-19   \>  .10e-22   \>  .25e-26   \\
\  \ 2.5   \> \>  .21e-12   \>  .52e-16   \>  .87e-21   \>  .21e-24   \>  .52e-28   \\
\  \ 5.0   \> \>  .12e-16   \>  .28e-20   \>  .48e-25   \>  .12e-28   \>  .28e-32   \\
\   10.0   \> \>  .20e-20   \>  .48e-24   \>  .81e-29   \>  .20e-32   \>  .48e-36   \\
\   25.0   \> \>  .30e-25   \>  .73e-29   \>  .12e-33   \>  .30e-37   \>   .73e-41   \\
\end{tabbing}
}$$

\vspace{-1.4cm}
\begin{quote}
{\bf Table 5.2.}\quad Relative accuracy 
$\Delta(\mu,t)$ defined in \parref{n2} for the asymptotic series of Section 2.1.
\end{quote}

Finally we consider the expansions of Sections 2.4 and 2.5.
Let the left-hand 
side of (2.35) be denoted by $W(\mu,t)$. Then we define as the error in the
expansions
\begin{equation}
\Delta(\mu,t):=\left|{\sfrac12W(\mu,t)}-1\right|\label{n3}
\end{equation}
When we use the series in (2.29), (2.33) and (2.34) with 5 terms, we obtain
the results of Table 5.3. We again observe that the accuracy improves as $\mu$ or $t$
increase. This shows the double asymptotic poperty of the modified 
expansions of Sections 2.4 and 2.5.
$$\vbox{
\begin{tabbing}
 10000\=1000\=10000000\=10000000\=10000000\=10000000\=10000000\kill
           \>  $\mu$   \>\ \ \ \ 5 \>\ \ \ 10\ \>\  \ \ 25 \>\ \  \ 50 \> \ \ 100   \\
\ $\ t$    \> \>            \>            \>            \>            \>           \\
\  \ .00   \> \>  .32e-09   \>  .78e-13   \>  .13e-17   \>  .32e-21   \>  .78e-25   \\
\  \ .25   \> \>  .12e-09   \>  .28e-13   \>  .47e-18   \>  .12e-21   \>  .28e-25   \\
\  \ .50   \> \>  .45e-11   \>  .11e-14   \>  .19e-19   \>  .46e-23   \>  .11e-26   \\
\  \ .75   \> \>  .57e-11   \>  .14e-14   \>  .24e-19   \>  .58e-23   \>  .14e-26   \\
\  \ 1.0   \> \>  .27e-11   \>  .65e-15   \>  .11e-19   \>  .27e-23   \>  .65e-27   \\
\  \ 1.5   \> \>  .29e-13   \>  .70e-17   \>  .12e-21   \>  .29e-25   \>  .70e-29   \\
\  \ 2.0   \> \>  .20e-13   \>  .48e-17   \>  .81e-22   \>  .20e-25   \>  .48e-29   \\
\ \  2.5   \> \>  .43e-14   \>  .11e-17   \>  .18e-22   \>  .43e-26   \>  .11e-29   \\
\ \  5.0   \> \>  .45e-17   \>  .11e-20   \>  .18e-25   \>  .45e-29   \>  .11e-32   \\
\   10.0   \> \>  .16e-20   \>  .38e-24   \>  .64e-29   \>  .16e-32   \>  .38e-36   \\
\end{tabbing}
}$$

\vspace{-1.4cm}
\begin{quote}
{\bf Table 5.3.}\quad Relative accuracy 
$\Delta(\mu,t)$ defined in \parref{n3} for the asymptotic series of Sections 2.4 
and 2.5.
\end{quote}

\section{Concluding remarks}%
As mentioned in Section 1.1 of the Introduction, several sources for
numerical algorithms for evaluating parabolic cylinder functions
are available in the literature, but not so many algorithms
make use of asymptotic expansions. The paper \cite{olpar}
is a rich source for asymptotic expansions, for all combinations of 
real and complex parameters, where always $|a|$ has to be large.
There are no published algorithms that make use of Olver's
expansions, although very efficient algorithms can be designed 
by using the variety of these expansions; \cite{blannew}
 is the only reference we found in which Olver's 
expansions are used for numerical computations.

We started our efforts in making algorithms 
for the case of real parameters. We selected 
appropriate expansions from
Olver's paper and for some cases we modified 
Olver's expansions in order to get expansions
having a double asymptotic property. A serious point
is making efficient use of the powerful Airy-type expansions that
are valid near the turning points of the differential equation
(and in much larger intervals and domains of the complex plane).
In particular, constructing reliable software for all possible 
combinations of the complex parameters $a$ and $z$ is a challengeing 
problem.

A point of research interest is also the construction of error bounds for
Olver's expansions and the modified expansions. Olver's paper is 
written before he developed the construction of bounds for the
remainders, which he based on methods for differential equations, 
and which are available now in his book \cite{olver}.
\section{Appendix: Computing the coefficients $f_k(\lambda)$ of (4.13)}%
We give the details on the computation of the coefficients $f_k(\lambda)$ that are used
in (4.13). The first step is to obtain coefficients $d_k$ in the expansion
\begin{equation}
w=d_0+d_1(t-\lambda)+d_2(t-\lambda)^2+\ldots,\label{x1}
\end{equation}
where $d_0=w_0$. From (4.6) we obtain
\begin{equation}
\frac{dw}{dt}=\frac wt\,\frac{t-\lambda}{w^2+w-\lambda}.\label{x2}
\end{equation}
Substituting \parref{x1} we obtain
\begin{equation}
d_1^2=\frac{w_0}{\lambda(1+2w_0)},\label{x3}
\end{equation}
where the saddle point $w_0$ is defined in (4.5). From the conditions
on the mapping (4.6) it follows that $d_1>0$. Higher order coefficients
$d_k$ can be obtained from the first ones by recursion.

When we have determined the coefficients in \parref{x1}
we can use (4.8) to obtain the coefficients $a_n(\lambda)$ of (4.10).
%
%
%

For computing in this way a set of coefficients $f_k(\lambda)$, say
$f_0(\lambda),\ldots,f_{15}(\lambda)$, we need more than 35 coefficients $d_k$
in \parref{x1}. Just taking
the square root in \parref{x3} gives for higher coefficients $d_k$ very
complicated expressions, and 
even by using computer algebra programs, as Maple, we
need suitable methods in computing the coefficients.

The computation of the coefficients
$d_k, a_n(\lambda)$ and $f_k(\lambda)$ is done
with a new parameter $\theta\in[0,\sfrac12\pi)$ which is defined by
\begin{equation}
4\lambda=\tan^2\theta.\label{x4}
\end{equation}
We also write
\begin{equation}
\sigma=\cos^2\sfrac12\theta, \label{x5}
\end{equation}
which is introduced earlier in (4.20) and (4.23). Then
\begin{equation}
w_0=\frac{1-\sigma}{2\sigma-1},\quad \lambda=\frac{\sigma(1-\sigma)}{(2\sigma-1)^2},\quad
d_1=\frac{2\sigma-1}{\sigma}.\label{x6}
\end{equation}
In particular the expressions for $w_0$ and $d_1$ are quite convenient,
because we can proceed without square roots in the computations.
Higher coefficients $d_k$ can be obtained by using \parref{x2}. 

The first relation $f_0(\lambda)=a_0^{(0)}(\lambda)=1$ easily follows from (4.3),
(4.8), \parref{x1} and \parref{x6}:
$$f_0(\lambda)={(1+4\lambda)^{\sfrac14}}\,\sqrt{{\frac {\lambda}{ w_0}
}}\,d_1=1.
$$

Then using (4.8) we obtain
$$a_0(\lambda)=1,\ a_1(\lambda)= -\frac{\cos^2\theta(1+2c)^2}{6(c+1)c^2},\
a_2(\lambda)=\frac{\cos^4\theta(20c^4+40c^3+30c^2+12c+3)}{24(c+1)^2c^4},$$
where $c=\sqrt{{\sigma}}=\cos\sfrac12\theta$. 
Using the scheme leading to (4.17)  one obtains the
coefficients $f_k(\lambda)$.  The first few coefficients are given in (4.19).

We observe that $f_k(\lambda)$ is a polynomial of degree $2k$ in $\sigma$
multiplied with $\rho^k$. If $a$ and $z$ are positive then $\sigma\in[0,1]$.
It follows that the sequence
$\{\rho^k/z^{2k}\}$ is an asymptotic scale when one or both
parameters $a$ and $z$ are large, and, hence, that $\{f_k(\lambda)/z^{2k}\}$  of
(4.13) is an asymptotic scale when one or both
parameters $a$ and $z$ are large.

Because of the relation in (4.27) and
$\widetilde\phi_k(\sigma)=\phi_k(\widetilde\tau)$, higher coefficients $f_k(\lambda)$ can also
be obtained from the recursion relation (2.11), which is obtained by using the
differential equation of the parabolic cylinder functions.

\end{document}